\documentclass[12pt,a4paper]{article}
\usepackage[T1]{fontenc}
\usepackage{mathptmx}
\usepackage{courier}
\usepackage[dvips]{graphicx}
\usepackage{psfrag}
\usepackage{amsmath,amssymb}        
\usepackage{latexsym}
\usepackage{url}
\usepackage[light]{kpfonts}

\def\ds{\displaystyle}

\title{A Non-Iterative Transformation Method\\ for Boundary-Layer with Power-Law Viscosity for Newtonian Fluids}
\author{Riccardo Fazio\\
Department of Mathematics, Computer Science,\\ Physical Sciences and Earth Sciences,\\
University of Messina \\
Viale F. Stagno D'Alcontres, 31 \\
98166 Messina, Italy \\
E-mails: \url{rfazio@unime.it} \\
Home-page: \url {http://mat521.unime.it/fazio}} 
\date{\today}
\begin{document}
\maketitle
\begin{abstract}
In this paper, we have defined and applied a non-ITM to an extended Blasius problem describing a 2D laminar boundary-layer with power-law viscosity for Newtonian fluids.
For a particular value of the parameter involved, this problem reduces to the celebrated Blasius problem and in this particular case, our method reduces to the T\"opfer non-iterative algorithm. 
In this case, we are able to compare favourably the obtained numerical result for the so-called missing initial condition with those available in the literature.
Moreover, we have listed the computed values of the missing initial condition for a large range of the parameter involved, and for illustrative purposes, we have plotted, for two values of the related parameter, the numerical solution computed rescaling a reference solution. 
\end{abstract}
\bigskip

\bigskip


\section{Introduction}
In this study, we are going to consider an extended version of the classical Blasius problem of boundary layer theory.
As well known, the basis of boundary layer theory was given by Prandtl at the beginning of the last century, see \cite{Prandtl:1904:UFK}.
In this context, T\"opfer \cite{Topfer:1912:BAB} in 1912 published a paper where he reduced the solution of the Blasius problem to the solution of a related initial value problem (IVP). 
This was the first definition and application of a non-iterative transformation method (ITM). 
As a consequence, non-ITMs have been applied to several problems of practical interest within the applied sciences. 
In fact, a non-ITM was applied to the Blasius equation with slip boundary condition, arising within the study of gas and liquid flows at the micro-scale regime \cite{Gad-el-Hak:1999:FMM,Martin:2001:BBL}, see \cite{Fazio:2009:NTM}.
Moreover, a non-ITM was applied also to the Blasius equation with moving wall considered by Ishak et al. \cite{Ishak:2007:BLM} or surface gasification studied by Emmons \cite{Emmons:1956:FCL} and recently by Lu and Law \cite{Lu:2014:ISB} or slip boundary conditions investigated by Gad-el-Hak \cite{Gad-el-Hak:1999:FMM} or Martin and Boyd \cite{Martin:2001:BBL}, see Fazio \cite{Fazio:2016:NIT} for details.
In particular, within these applications, we found a way to solve non-iteratively the Sakiadis problem \cite{Sakiadis:1961:BLBa,Sakiadis:1961:BLBb}.
The application of a non-ITM to an extended Blasius problem has been the subject of a recent manuscript \cite{Fazio:2020:NIT}.
The interested reader can find in \cite{Fazio:2019:NIT} a recent review dealing with the non-ITM and its applications.

The non-ITM is based on scaling invariance theory.
For its application, the governing differential equation as well as the prescribed initial conditions have to be invariant with respect to a scaling group of point transformations.
Of course, there are several problems in the applied sciences that luck this kind of invariance and consequently cannot be solved by a non-ITM.
However, in all those cases we can use an iterative extension of our method.
In fact, an iterative extension of T\"opfer's algorithm has been introduced, for the numerical solution of free BVPs, by Fazio \cite{Fazio:1990:SNA}. 
This iterative extension has been applied to several problems of interest: free boundary problems \cite{Fazio:1990:SNA,Fazio:1997:NTE,Fazio:1998:SAN},
a moving boundary hyperbolic problem \cite{Fazio:1992:MBH}, Homann and Hiemenz problems governed by the Falkner-Skan equation and a mathematical model describing the study of the flow of an incompressible fluid around a slender parabola of revolution \cite{Fazio:1994:FSE,Fazio:1996:NAN},
one-dimensional parabolic moving boundary problems \cite{Fazio:2001:ITM}, two variants of the Blasius problem \cite{Fazio:2009:NTM}, namely: a boundary layer problem over moving surfaces, studied first by Klemp and Acrivos \cite{Klemp:1972:MBL}, and a boundary layer problem with slip boundary condition, that has found application in the study of gas and liquid flows at the micro-scale regime \cite{Gad-el-Hak:1999:FMM,Martin:2001:BBL}, parabolic problems on unbounded domains \cite{Fazio:2010:MBF} and, recently, see \cite{Fazio:2015:ITM}, a further variant of the Blasius problem in boundary layer theory: the so-called Sakiadis problem \cite{Sakiadis:1961:BLBa,Sakiadis:1961:BLBb}.
Moreover, this iterative extension can be used to investigate the existence and uniqueness question for different classes of problems, as shown for free BVPs in \cite{Fazio:1997:NTE}, and for problems in boundary layer theory in \cite{Fazio:2020:EUB}.
A recent review dealing with, the derivation and applications of the ITM can be found, by the interested reader, in \cite{Fazio:2019:ITM}.
A unifying framework, providing proof that the non-ITM is a special instance of the ITM and consequently can be derived from it, has been the argument of the paper \cite{Fazio:2020:SIT}.

The mathematical model arising in the study of a 2D laminar boundary-layer with power-law viscosity for Newtonian fluids is given by, see Schlichting and Gersten \cite{Schlichting:2000:BLT} or Benlahsen et al. \cite{Benlahsen:2008:GBE},
\begin{align}\label{eq:goveq} 
& {\displaystyle \frac{d}{d\eta}\left(\left|\frac{d^{2}f}{d\eta^2}\right|^{n-1} \frac{d^{2}f}{d\eta^2}\right)} + \frac{1}{n+1} f {\displaystyle \frac{d^{2}f}{d\eta^2}} = 0 \nonumber \\[-1ex]
&\\[-1ex]
& f(0) = {\displaystyle \frac{df}{d\eta}}(0) = 0 \ , \qquad {\displaystyle \frac{df}{d\eta}}(\eta) \rightarrow 1 \quad \mbox{as} \quad \eta \rightarrow \infty \ , \nonumber 
\end{align}
where $f(\eta)$ is the non-dimensional stream function and $n$ is a given positive value bigger than zero.
Let us remark here, that when  $n = 1$ the BVP (\ref{eq:goveq}) reduces to the celebrated Blasius problem \cite{Blasius:1908:GFK}.

\section{The non-ITM} 
In order to define our non-ITM we need to require the invariance of the governing differential equation and of the prescribed initial conditions in (\ref{eq:goveq}) with respect to the scaling group of point transformation
\begin{equation}\label{eq:scaling2}
f^* = \lambda f \ , \qquad \eta^* = \lambda^{\delta} \eta \ .   
\end{equation}
It is easily seen, that the governing differential equation and the prescribed initial conditions are invariant on the condition that $\delta = (2-n)/(1-2n)$.
Therefore, we have to require that $n \ne 1/2$ and $n \ne 2$.
Now, we can integrate the governing equation in (\ref{eq:goveq}) written in the star variables on $[0, \eta^*_\infty]$, where $\eta^*_\infty$ is a suitable truncated boundary, with initial conditions
\begin{equation}\label{eq:ICs2}
f^*(0) = \frac{df^*}{d\eta^*}(0) = 0 \ , \quad {\displaystyle \frac{d^2f^*}{d\eta^{*2}}(0)} = 1 \ ,
\end{equation}
in order to compute an approximation $\frac{df^*}{d\eta^*}(\eta^*_{\infty})$ for $\frac{df^*}{d\eta^*}(\infty)$ and the corresponding value of $\lambda$ by the equation
\begin{equation}\label{eq:lambda}
\lambda = \left[\frac{df^*}{d\eta^*}(\eta_\infty^*)\right]^{1/(1-\delta)} \ .
\end{equation}
Once the value of $\lambda$ has been computed by equation (\ref{eq:lambda}), we can find the missed initial condition by the relation
\begin{equation}\label{eq:MIC}
\frac{d^2f}{d\eta^{2}}(0) =  \lambda^{2\delta-1} \frac{d^2f^*}{d\eta^{*2}}(0) \ ,
\end{equation}
and rescale the solution components according to
\begin{equation}\label{eq:Rescale}
f(\eta) = \lambda^{-1} f^*(\eta^*)\ , \quad \frac{df}{d\eta}(\eta) =  \lambda^{\delta-1} \frac{df^*}{d\eta^*}(\eta^*) \ , \quad \frac{d^2f}{d\eta^{2}}(\eta) =  \lambda^{2\delta-1} \frac{d^2f^*}{d\eta^{*2}}(\eta^*) \ .
\end{equation}

In the non-ITM we proceed as follows: we set the values of $\eta_{\infty}^*$ and integrate the IVP, the governing differential equation and the initial conditions (\ref{eq:ICs2}), on $[0, \eta_{\infty}^*]$.

\section{Numerical results}
In this section, after some sample computations we set the truncated boundary value equal to ten, that is for all values of $n$ we use $\eta^*_\infty = 10$.
As mentioned before, the case $n = 1$ is the Blasius problem and in this case, our non-ITM reduces to the classical T\"opfer algorithm, see T\"opfer \cite{Topfer:1912:BAB}.
In this case, the numerical value computed for the missing initial condition, namely $0.332057268052$, can be compared with the value $0.332057336215$ obtained by Fazio \cite{Fazio:1992:BPF} by a free boundary formulation of the Blasius problem or the value $0.33205733621519630$ computed by Boyd \cite{Boyd:1999:BFC} who believes that all the decimal digits are correct.

\begin{figure}[!hbt]
	\centering
\psfrag{e}[l][]{$\eta$,$\eta^*$}  
\psfrag{df}[l][]{$\ds \frac{df}{d\eta}$} 
\psfrag{ddf}[l][]{$\ds \frac{d^2f}{d\eta^2}$} 
\psfrag{df*}[l][]{$\ds \frac{df^*}{d\eta^*}$} 
\psfrag{ddf*}[l][]{$\ds \frac{d^2f^*}{d\eta^{2*}}$} 
\includegraphics[width=12cm,height=8cm]{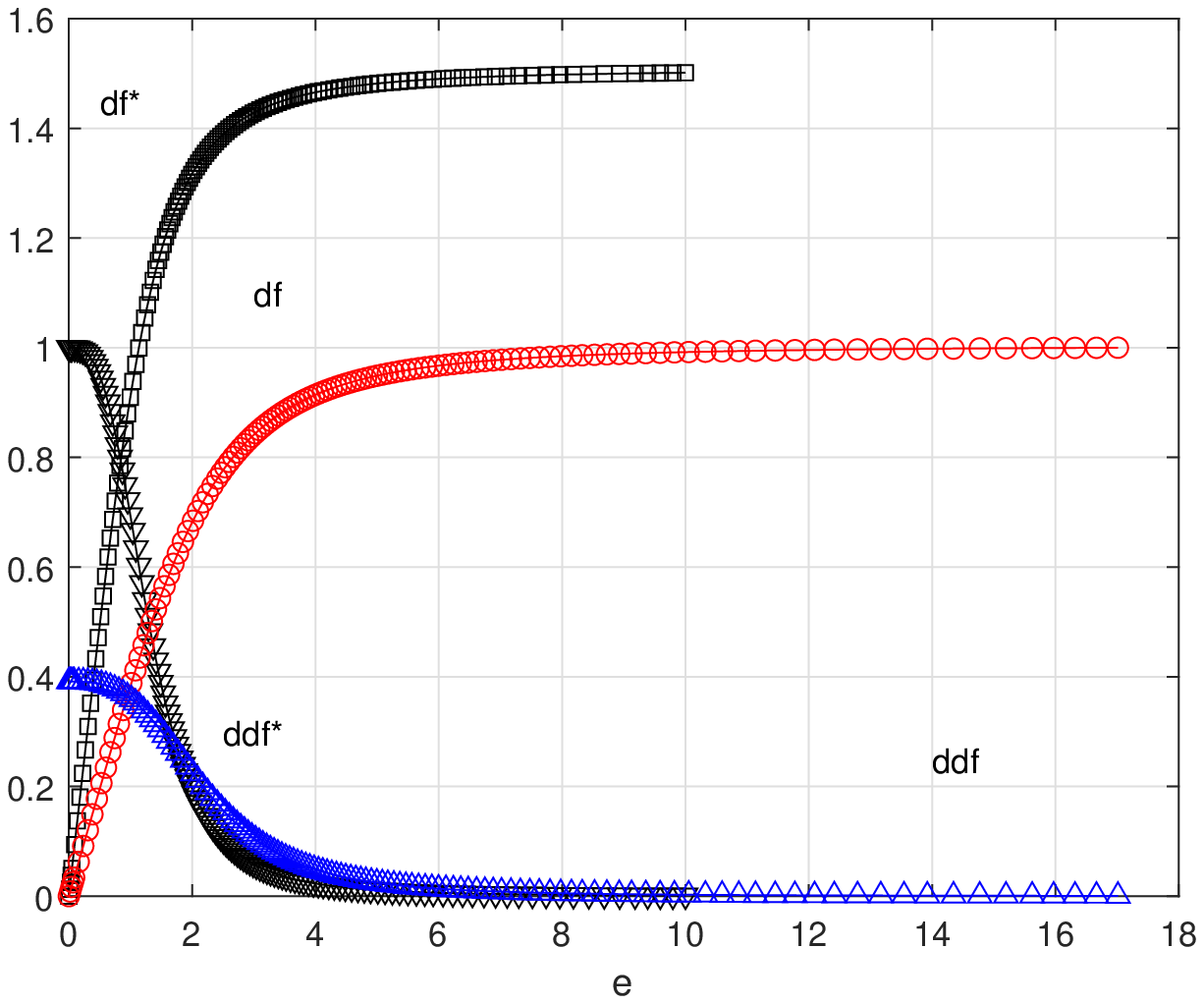} \\[1ex]
\includegraphics[width=12cm,height=8cm]{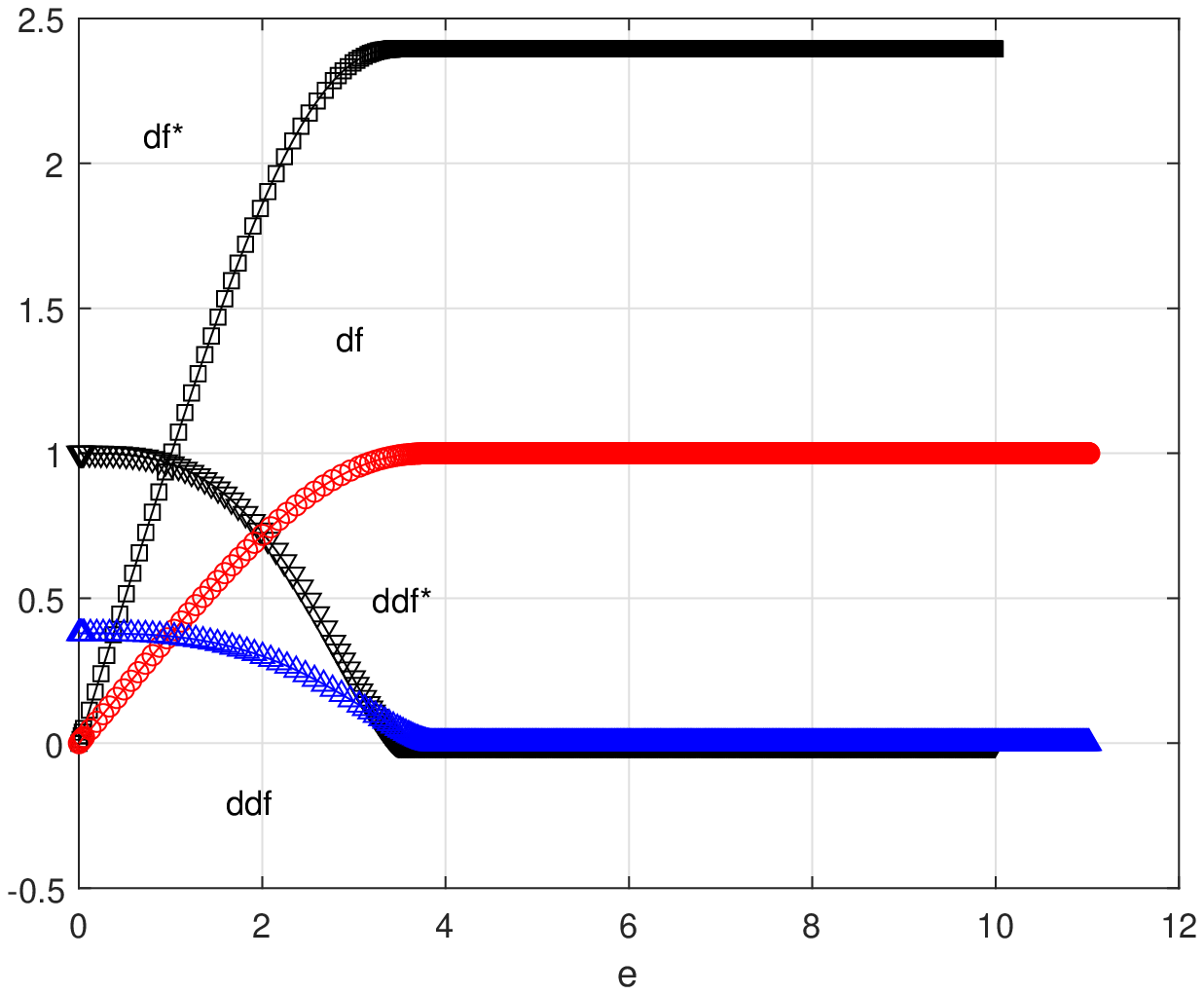}
\caption{Numerical results of the non-ITM for (\ref{eq:goveq}) with: top frame $n = 0.3$ and bottom frame $n= 1.7$.}
	\label{fig:1n03n17}
\end{figure}

In table \ref{tab:NITM:missingIC} we report the chosen parameter values and the related missing initial conditions $\frac{d^2f}{d\eta^2}(0)$.
\begin{table}[!hbt]
\caption{Numerical data and results.}
\vspace{.5cm}
\renewcommand\arraystretch{1.3}
	\centering
		\begin{tabular}{lr@{.}llr@{.}l}
\hline \\[-2.2ex]
{$n$} & \multicolumn{2}{c}%
{$ {\displaystyle \frac{d^2f}{d\eta^2}(0)}$} 
& {$n$} & \multicolumn{2}{c}%
{$ {\displaystyle \frac{d^2f}{d\eta^2}(0)}$}\\[1.5ex]
\hline
0.1       & 0 & 826474 & 1.1 & 0 & 337833 \\ 
0.2       & 0 & 490340 & 1.2 & 0 & 344165 \\ 
0.3       & 0 & 391514 & 1.3 & 0 & 350851 \\ 
0.4       & 0 & 350395 & 1.4 & 0 & 357752 \\ 
0.5       & 0 & 337170 & 1.5 & 0 & 364772 \\ 
0.6       & 0 & 323945 & 1.6 & 0 & 371841 \\ 
0.7	    	& 0 & 322033 & 1.7 & 0 & 378906 \\ 
0.8       & 0 & 323543 & 1.8 & 0 & 385936 \\ 
0.9	    	& 0 & 327139 & 1.9 & 0 & 392894 \\ 
1         & 0 & 332057 & 2   & 0 & 399852 \\ 
\hline			
		\end{tabular}
	\label{tab:NITM:missingIC}
\end{table}
The values listed in table \ref{tab:NITM:missingIC} for the two values of $n = 1/2$ and $n = 2$ are second order approximations.
To verify the last value reported in this table we have also considered the case $n = 1.999$ and found the missing initial condition value $0.399700$.
Figure \ref{fig:1n03n17} shows the solution of the extended Blasius problem in the particular case when we set $n = 0.3$ and $n = 1.7$.

\section{Concluding remarks} 
In this paper we have defined and applied a non-ITM to an extended Blasius problem describing a 2D laminar boundary-layer with power-law viscosity for Newtonian fluids as described by Schlichting and Gersten \cite{Schlichting:2000:BLT} or Benlahsen et al. \cite{Benlahsen:2008:GBE}.
This problem reduces to the celebrated Blasius problem for $n = 1$ and in this particular case, our method reduces to the T\"opfer non-iterative algorithm \cite{Topfer:1912:BAB}. 
In this case, we are able to compare favourably the obtained numerical result for the so-called missing initial condition with those available in the literature.
Moreover, we have listed the computed values of the missing initial condition for a large range of the parameter involved, and for illustrative purposes, we have plotted, for two values of the related parameter, the numerical solution computed rescaling a reference solution. 

\vspace{1.5cm}

\noindent {\bf Acknowledgement.} {The research of this work was 
partially supported by the FFABR grant of the University of Messina and by the GNCS of INDAM.}


\begin{thebibliography}{10}

\bibitem{Benlahsen:2008:GBE}
M.~Benlahsen, M.~Guedda, and R.~Kersner.
\newblock The generalized {B}lasius equation revisited.
\newblock {\em Math. Computer Model.}, 47:1063--1076, 2008.

\bibitem{Blasius:1908:GFK}
H.~Blasius.
\newblock Grenzschichten in {F}l\"{u}ssigkeiten mit kleiner {R}eibung.
\newblock {\em Z. Math. Phys.}, 56:1--37, 1908.

\bibitem{Boyd:1999:BFC}
J.~P. Boyd.
\newblock The {B}lasius function in the complex plane.
\newblock {\em Exp. Math.}, 8:381--394, 1999.

\bibitem{Gad-el-Hak:1999:FMM}
M.~Gad el~Hak.
\newblock The fluid mechanics of microdevices --- the {F}reeman scholar
  lecture.
\newblock {\em J. Fluids Eng.}, 121:5--33, 1999.

\bibitem{Emmons:1956:FCL}
H.~W. Emmons.
\newblock The film combustion of liquid fluid.
\newblock {\em ZAMM - J. Appl. Math. Mech.}, 36:60--71, 1956.

\bibitem{Fazio:1992:BPF}
R.~Fazio.
\newblock The {Blasius} problem formulated as a free boundary value problem.
\newblock {\em Acta Mech.}, 95:1--7, 1992.

\bibitem{Fazio:1992:MBH}
R.~Fazio.
\newblock A moving boundary hyperbolic problem for a stress impact in a bar of
  rate-type material.
\newblock {\em Wave Motion}, 16:299--305, 1992.

\bibitem{Fazio:1994:FSE}
R.~Fazio.
\newblock The {F}alkner-{S}kan equation: numerical solutions within group
  invariance theory.
\newblock {\em Calcolo}, 31:115--124, 1994.

\bibitem{Fazio:1996:NAN}
R.~Fazio.
\newblock A novel approach to the numerical solution of boundary value problems
  on infinite intervals.
\newblock {\em SIAM J. Numer. Anal.}, 33:1473--1483, 1996.

\bibitem{Fazio:1997:NTE}
R.~Fazio.
\newblock A numerical test for the existence and uniqueness of solution of free
  boundary problems.
\newblock {\em Appl. Anal.}, 66:89--100, 1997.

\bibitem{Fazio:1998:SAN}
R.~Fazio.
\newblock A similarity approach to the numerical solution of free boundary
  problems.
\newblock {\em SIAM Rev.}, 40:616--635, 1998.

\bibitem{Fazio:2001:ITM}
R.~Fazio.
\newblock The iterative transformation method: numerical solution of
  one-dimensional parabolic moving boundary problems.
\newblock {\em Int. J. Computer Math.}, 78:213--223, 2001.

\bibitem{Fazio:2009:NTM}
R.~Fazio.
\newblock Numerical transformation methods: {B}lasius problem and its variants.
\newblock {\em Appl. Math. Comput.}, 215:1513--1521, 2009.

\bibitem{Fazio:2015:ITM}
R.~Fazio.
\newblock The iterative transformation method for the {S}akiadis problem.
\newblock {\em Comput. \& Fluids}, 106:196--200, 2015.

\bibitem{Fazio:2016:NIT}
R.~Fazio.
\newblock A non-iterative transformation method for {B}lasius equation with
  moving wall or surface gasification.
\newblock {\em Int. J. Non-Linear Mech.}, 78:156--159, 2016.

\bibitem{Fazio:2019:ITM}
R.~Fazio.
\newblock The iterative transformation method.
\newblock {\em Int. J. Non-Linear Mech.}, 116:181--194, 2019.

\bibitem{Fazio:2019:NIT}
R.~Fazio.
\newblock The non-iterative transformation method.
\newblock {\em Int. J. Non-Linear Mech.}, 114:41--48, 2019.

\bibitem{Fazio:2020:EUB}
R.~Fazio.
\newblock Existence and uniqueness of bvps defined on infinite intervals:
  Insight from the iterative transformation method, 2020.
\newblock Preprint available at the URL:
  \url{http://mat521.unime.it/~fazio/preprints/NTest2020.pdf}.

\bibitem{Fazio:2020:NIT}
R.~Fazio.
\newblock A non-iterative transformation method for an extended {B}lasius
  problem, 2020.
\newblock Preprint available at the URL:
  \url{http://mat521.unime.it/~fazio/preprints/ExBlasius2020.pdf}.

\bibitem{Fazio:2020:SIT}
R.~Fazio.
\newblock Scaling invarance theory and numerical transformation methods: A
  unifying framework, 2020.
\newblock Preprint available at the URL:
  \url{http://mat521.unime.it/~fazio/preprints/SINTMUF2020.pdf}.

\bibitem{Fazio:1990:SNA}
R.~Fazio and D.~J. Evans.
\newblock Similarity and numerical analysis for free boundary value problems.
\newblock {\em Int. J. Computer Math.}, 31:215--220, 1990.
\newblock 39~:~249, 1991.

\bibitem{Fazio:2010:MBF}
R.~Fazio and S.~Iacono.
\newblock On the moving boundary formulation for parabolic problems on
  unbounded domains.
\newblock {\em Int. J. Computer Math.}, 87:186--198, 2010.

\bibitem{Ishak:2007:BLM}
A.~Ishak, R.~Nazar, and I.~Pop.
\newblock Boundary layer on a moving wall with suction and injection.
\newblock {\em Chin. Phys. Lett.}, 24:2274--2276, 2007.

\bibitem{Klemp:1972:MBL}
J.~P. Klemp and A.~Acrivos.
\newblock A moving-wall boundary layer with reverse flow.
\newblock {\em J. Fluid Mech.}, 53:177--191, 1972.

\bibitem{Lu:2014:ISB}
Z.~Lu and C.~K. Law.
\newblock An iterative solution of the {B}lasius flow with surface
  gasification.
\newblock {\em Int. J. Heat and Mass Transfer}, 69:223--229, 2014.

\bibitem{Martin:2001:BBL}
M.~J. {Martin} and I.~D. {Boyd}.
\newblock {Blasius boundary layer solution with slip flow conditions}.
\newblock In {\em Rarefied Gas Dynamics: 22nd International Symposium}, volume
  585 of {\em American Institute of Physics Conference Proceedings}, pages
  518--523, 2001, DOI: 10.1063/1.1407604.

\bibitem{Prandtl:1904:UFK}
L.~Prandtl.
\newblock {\"U}ber {F}l{\"u}ssigkeiten mit kleiner {R}eibung.
\newblock In {\em Proc. Third Inter. Math. Congr.}, pages 484--494, 1904.
\newblock Engl. transl. in {NACA} Tech. Memo. 452.

\bibitem{Sakiadis:1961:BLBa}
B.~C. Sakiadis.
\newblock Boundary-layer behaviour on continuous solid surfaces: I.
  {B}oundary-layer equations for two-dimensional and axisymmetric flow.
\newblock {\em AIChE J.}, 7:26--28, 1961.

\bibitem{Sakiadis:1961:BLBb}
B.~C. Sakiadis.
\newblock Boundary-layer behaviour on continuous solid surfaces: {II}. {T}he
  boundary layer on a continuous flat surface.
\newblock {\em AIChE J.}, 7:221--225, 1961.

\bibitem{Schlichting:2000:BLT}
H.~Schlichting and K.~Gersten.
\newblock {\em Boundary Layer Theory}.
\newblock Springer, Berlin, 8th edition, 2000.

\bibitem{Topfer:1912:BAB}
K.~T{\"o}pfer.
\newblock Bemerkung zu dem {A}ufsatz von {H}. {B}lasius: {G}renzschichten in
  {F}l{\"u}ssigkeiten mit kleiner {R}eibung.
\newblock {\em Z. Math. Phys.}, 60:397--398, 1912.

\end{thebibliography}
\end{document}